\documentclass[12pt,leqno]{article}
\usepackage{amssymb,eufrak}
\begin{document}
\title{On Jacobi Sums in 
$\mathbb Q(\zeta_p )$}

\author{
Bruno Angl\`  es \\
Universit\'e de Caen, 
 CNRS UMR 6139, \\
Campus II, Boulevard Mar\'echal Juin, 
B.P. 5186, \\
14032 Caen Cedex, France.\\
E-mail: bruno.angles@math.unicaen.fr
\and 
Filippo A. E. Nuccio\\
Istituto Guido Castelnuovo,
Universit\`a "La Sapienza",\\
5, Piazzale Aldo Moro,\\
00186 Rome, Italy.\\
E-mail: nuccio@mat.uniroma1.it}
\date{21 november 2007}
\maketitle
\begin{abstract}
We study the $p$-adic behavior of Jacobi sums for $\mathbb Q(\zeta_p)$ and link this study to the $p$-Sylow subgroup of the class group of $\mathbb Q(\zeta_p)^+$ and to  some properties of the jacobian of the Fermat curve $X^p+Y^p=1$ over $\mathbb F_{\ell}$ where $\ell$ is a prime number distinct from $p.$
 \end{abstract}\par
 Let $p$ be a prime number, $p\geq 5.$ Iwasawa has shown that the $p$-adic properties of Jacobi sums for $\mathbb Q(\zeta_p)$ are linked to Vandiver's Conjecture (see \cite{IWA1}). In this paper, we follow Iwasawa's ideas and study the $p$-adic properties of the subgroup $J$  of $\mathbb Q(\zeta_p)^*$ generated by Jacobi sums.\par 
 Let $A$ be the $p$-Sylow subgroup of the class group of $\mathbb Q(\zeta_p).$ If $E$ denotes the group of units of $\mathbb Q(\zeta_p),$ then if Vandiver's Conjecture is true for $p,$ by Kummer  Theory,  we must have $\frac{A^-}{pA^-}\hookrightarrow {\rm Gal}(\mathbb Q(\zeta_p)({}^p\sqrt{E})/\mathbb Q(\zeta_p)).$ Note that $J$ is analoguous for the odd part to the group of cyclotomic units for the even part. We introduce a submodule $W$ of $\mathbb Q(\zeta_p)^*$ which was already considered by Iwasawa (\cite{IWA2}). This module can be thought as the analogue  for the odd part of the group of units for the even  part.  We observe that $J\subset W$ and if the Iwasawa-Leopoldt Conjecture is true for $p$ then $W\mathbb (Q(\zeta_p)^*)^p=J(Q(\zeta_p)^*)^p.$ We prove that if $pA^-=\{ 0\}$ then (Corollary \ref{Corollary4}):
 $$\frac{A^+}{pA^+}\hookrightarrow {\rm Gal}(\mathbb Q(\zeta_p)({}^p\sqrt{W})/\mathbb Q(\zeta_p)).$$
 The last part of our paper is devoted to the study of the jacobian of the Fermat curve $X^p+Y^p=1$ over $\mathbb F_{\ell}$ where $\ell$ is a prime number, $\ell \not =p.$ It is well-known that Jacobi sums play an important role in the study of that jacobian. Following ideas developped by Greenberg (\cite{GRE}), we prove that Vandiver's Conjecture is equivalent to some properties  of that jacobian (for the precise statement see Corollary \ref{Corollary5}).\par
 The authors thank Cornelius Greither for interesting  discussions on  the converse of Kummer's Lemma which led us to the study of the analoguous statement for the odd part . The second author thanks the mathematicians of the Laboratoire de Math\'ematiques Nicolas Oresme for their hospitality during his stay at Caen. \par

 \section{Notations}
 ${}$\par
 Let $p$ be a prime number, $p\geq 5.$ Let $\zeta_{p}\in \mu_{p}\setminus \{ 1\},$ and let $L=\mathbb Q(\zeta_{p}).$ Set ${\cal O}=\mathbb Z[\zeta_{p}]$ and $E={\cal O}^{*}.$
 Let $\Delta ={\rm Gal}(L/\mathbb Q)$ and let $\widehat{\Delta}={\rm Hom}(\Delta , \mathbb Z_{p}^{*}).$ Let $\cal I$ be the group of fractional ideals of $L$ which are prime to $p,$ and let 
 $\cal P$ be the group of principal ideals in $\cal I.$ Let $A$ be the $p$-Sylow subgroup of the ideal class group of $L.$\par
 Set $\pi=\zeta_{p}-1,$ $K=\mathbb Q_{p}(\zeta_{p}),$ $U=1+\pi^{2}\mathbb Z_{p}[\zeta_{p}].$ Observe that if ${\cal A} \in {\cal P},$ then there exists $\alpha \in L^{*}\cap U$ such that ${\cal A}
 = \alpha {\cal O}.$ If $H$ is a subgroup of $U,$ we will denote the  closure of $H$ in $U$ by $\overline{H}.$ Let $\omega \in {\widehat \Delta}$ be the Teichm\" uller character, i.e. :
 $$\forall \sigma \in \Delta, \, \sigma (\zeta_{p})=\zeta_{p}^{\omega (\sigma)}.$$
 For $\rho \in {\widehat \Delta},$ we set:
 $$ e_{\rho}=\frac{1}{p-1}\sum_{\delta \in \Delta}\rho^{-1}(\delta) \delta \, \in \mathbb Z_{p}[\Delta].$$
 If $M$ is a $\mathbb Z_{p}[\Delta]$-module, for $\rho \in {\widehat \Delta},$ we set:
 $$ M(\rho )=e_{\rho}M.$$\par
 For $\psi \in {\widehat \Delta},$ $\psi$ odd, recall that:
 $$ B_{1,\psi}=\frac{1}{p}\sum_{a=1}^{p-1}a\psi (a).$$
 Set:
 $$ \theta=\frac{1}{p}\sum_{a=1}^{p-1}a\sigma_{a}^{-1}\, \in \mathbb Q [\Delta],$$
 where $\sigma_{a} \in \Delta $ is such that $\sigma_{a}(\zeta_{p})=\zeta_{p}^{a}.$ Observe that we have the following equality in $\mathbb C [\Delta] :$
 $$ \theta = \frac{N}{2}+\sum_{\psi \in {\widehat \Delta},\, \psi \, {\rm odd}} B_{1,\psi^{-1}} e_{\psi},$$
 where $N=\sum_{\delta \in \Delta} \delta .$\par
 Let $M$ be a $\mathbb Z[\Delta]$-module, we set:
 $$ M^{-}=\{ m\in M,\, \sigma_{-1}(m)=-m \},$$
 $$ M^{+}= \{ m\in M,\, \sigma_{-1}(m)=m \}.$$
 If $M$ is an abelian group of finite type, we set:
 $$ M[p]=\{ m\in M,\, pm=0\},$$
 $$ d_{p}M={\rm dim}_{\mathbb F_{p}} \frac{M}{pM}\,.$$
 
 \section{Background on Jacobi Sums}
 ${}$\par
Let $ Cl (L)$ be the ideal class group of $L,$ then:
$$ Cl(L)\simeq \frac{{\cal I}}{{\cal P}}.$$
Note that we have a natural $\mathbb Z[\Delta]$-morphism (see \cite{IWA2}, pages 102-103):
$$\phi: ({\rm Ann}_{\mathbb Z[\Delta]} Cl(L))^{-}\rightarrow {\rm Hom}_{\mathbb Z[\Delta]}(Cl(L), \frac{E^{+}}{(E^{+})^ 2}).$$
For the convenience of the reader, we recall the construction of $\phi.$ Let $x\in ({\rm Ann}_{\mathbb Z[\Delta]} Cl(L))^{-}.$ Let ${\cal A}\in {\cal I},$ we have:
$$ {\cal A}^x= \gamma_{a}{\cal O},$$
where $\gamma_{a}\in L^*\cap U.$ Now:
$$ \overline{\gamma_{a}}=\varepsilon_{a}\gamma_{a}^{-1},$$
for some $\varepsilon_{a}\in E^{+}\cap U.$ One can prove that we obtain a well-defined morphism of $\mathbb Z[\Delta]$-modules:
$\phi(x): Cl(L)\rightarrow \frac {E^{+}} {(E^{+})^2},$ ${\rm class\, of \, } {\cal A}\mapsto {\rm class\, of\, } \varepsilon_{a}.$\par
\noindent In this paragraph, we will study the kernel of the morphism $\phi .$\par
Let $\cal W$ be the set of elements $f\in {\rm Hom}_{\mathbb Z[\Delta ]}({\cal I}, L^*)$ such that:\par
\noindent - $f({\cal I})\subset U,$\par
\noindent - there exists $\beta (f)\in \mathbb Z[\Delta ]$ such that for all $\alpha \in L^*\cap U,$ $f(\alpha {\cal O})= \alpha ^{\beta (f)}.$\par
\noindent One can prove that if $f\in {\cal W}$ then $\beta (f)$ is unique, the map $\beta : {\cal W}\rightarrow \mathbb Z[\Delta ]$ is an injective $\mathbb Z[\Delta ]$-morphism and $\beta ({\cal W})\subset {\rm Ann}_{\mathbb Z[\Delta ]}(Cl (L))$ (see \cite{ANG&BEL}). If ${\cal B}$ denotes the group of Hecke characters of type $(A_0)$  that have values in $\mathbb Q(\zeta_p)$ (see \cite{IWA2}), then one can prove that ${\cal B}$ is isomorphic to ${\cal W}.$\par
\newtheorem{Lemma1}{Lemma}[section]
 \begin{Lemma1} \label{Lemma1}
${\rm Ker}\, \phi = \beta ({\cal W}^{-}).$
\end{Lemma1}
\noindent {\sl Proof}  We just prove the inclusion ${\rm Ker}\, \phi \subset  \beta ({\cal W}^{-}).$ Let $x \in {\rm Ker}\, \phi .$ Let ${\cal A}\in {\cal I},$ then there exists an unique $\gamma_{a}\in L^*\cap U$ such that $\overline{\gamma_{a}}\gamma_{a}=1$ and:
$$ {\cal A}^x=\gamma_{a} {\cal O}.$$
Let $f:{\cal I}\rightarrow L^*,$ ${\cal A}\mapsto \gamma_{a}.$ It is not difficult to see that $f\in {\rm Hom}_{\mathbb Z[\Delta ]}({\cal I}, L^*)$ and $f({\cal I})\subset U.$ Now, let $\alpha \in L^*\cap U,$ we have:
$$ f(\alpha {\cal O})= \alpha ^{x} u,$$
for some $u\in E.$ Since $x\in \mathbb Z[\Delta ]^{-}$ and $\alpha, f(\alpha {\cal O})\in U,$ we must have $u=1.$ Threfore $f\in {\cal W}^{-}$ and $x=\beta (f).$ $\diamondsuit$\par
Now, we recall some basic properties of Gauss and Jacobi sums (we refer the reader to \cite{WAS}, paragraph 6.1).\par
\noindent Let $P$ be a prime ideal in ${\cal I}$ and let $\ell$ be the prime number such that ${\ell}\in P.$ We fix $\zeta_{\ell }\in \mu_{\ell }\setminus \{ 1\}.$
 Set $\mathbb F_P=\frac{{\cal O}}{P}.$ Let $\chi_P: \mathbb F_P^*\rightarrow \mu_p,$ such that:
$$\forall \alpha \in \mathbb F_P^*,\, \chi_P(\alpha )\equiv \alpha ^{\frac{1-NP}{p}}\pmod{P},$$
where $NP=\mid \frac {{\cal O}}{P}\mid .$ For $a\in \frac{\mathbb Z}{p\mathbb Z},$ we set:
$$\tau_a(P)=-\sum_{\alpha \in \mathbb F_P}\chi_P^{a}(\alpha )\zeta_{\ell }^{{\rm Tr}_{\mathbb F_P/\mathbb F_{\ell}}(\alpha )}.$$
We also set $\tau (P)=\tau_1(P).$ For $a,b\in \frac{\mathbb Z}{p\mathbb Z},$ we set:
$$j_{a,b}(P)=-\sum_{\alpha \in \mathbb F_P}\chi_P^{a}(\alpha )\chi_P^{b}(1-\alpha).$$
Then:\par
\noindent - if $a+b\equiv 0\pmod{p},$\par
i) if $a\not \equiv 0\pmod{p},$ $j_{a,b}(P)=1,$\par
ii) if $a\equiv 0\pmod{p},$ $j_{a,b}(P)=2-NP,$\par
\noindent - if $a+b\not \equiv 0\pmod{p},$ we have:
$$ j_{a,b}(P)=\frac{\tau_a (P)\tau_b(P)}{\tau_{a+b}(P)}.$$
Observe that $\tau (P) \equiv 1\pmod{\pi},$ and therefore (see \cite{IWA1}, Theorem 1):
$$\forall a,b \in \frac{\mathbb Z}{p\mathbb Z},\, j_{a,b}(P)\in U.$$
Let $\Omega$ be the compositum of the fields $\mathbb Q(\zeta_{\ell })$ where $\ell $ runs through the prime numbers distinct from $p.$ The map $P\mapsto \tau (P)$ induces by linearity a $\mathbb Z[\Delta ]$-morphism:
$$\tau : {\cal I}\rightarrow \Omega (\zeta_p )^*.$$
Let $\cal G$ be the sub-$\mathbb Z[\Delta ]$-module of ${\rm Hom}_{\mathbb Z[\Delta ]} ({\cal I}, \Omega (\zeta_p)^*)$ generated by $\tau.$ We set:
$${\cal J}={\cal G}\cap {\rm Hom}_{\mathbb Z[\Delta]} ({\cal I}, L^*).$$
Let $\cal S$ be the Stickelberger ideal of $L,$ i.e. : ${\cal S}=\mathbb Z[\Delta ]\theta \cap \mathbb Z[\Delta ].$ Then one can prove the following facts (see \cite{ANG&BEL}):\par
\noindent - ${\cal J}\subset {\cal W},$\par
\noindent - the map $\beta : {\cal W} \rightarrow \mathbb Z[\Delta ]$ induces an isomorphism of $\mathbb Z[\Delta ]$-modules :
$${\cal J}\simeq {\cal S}.$$
\newtheorem{Lemma2}[Lemma1]{Lemma}
 \begin{Lemma2} \label{Lemma2}
 Let ${\cal N}\in {\rm Hom}_{\mathbb Z[\Delta ]}(I_L, L^*)$ be the ideal norm map. Then, as a $\mathbb Z$-module:
 $${\cal J}={\cal N}\mathbb Z\oplus \bigoplus_{n=1}^{(p-1)/2}j_{1,n}\mathbb Z.$$
 \end{Lemma2}
 \noindent{\sl Proof} Recall that, for $1\leq n\leq  p-2,$ for a prime $P$ in ${\cal I},$  we have:
 $$j_{1,n}(P)=-\sum_{\alpha \in \mathbb F_P}\chi_P(\alpha )\chi_P^{n}(1-\alpha )=\frac{\tau (P)\tau_n(P)}{\tau_{n+1}(P)}.$$
 Thus, for $1\leq n\leq p-2,$ we have:
 $$j_{1,n}=\tau^{1+\sigma_n-\sigma _{1+n}}=\frac{\tau \tau_n}{\tau_{n+1}},$$
 where  for $a \in \mathbb F_p^*,$ $\tau^{\sigma_a}=\tau_a.$ Observe that:
 $$\forall a \in \mathbb F_p^*,\, \tau_a \tau_{-a}= N.$$
 Thus ${\cal N} \in {\cal J}.$ Since ${\cal J}\simeq {\cal S},$ ${\cal J}$ is a $\mathbb Z$-module of rank $(p+1)/2.$ It is not difficult to show that (see \cite{IWA1}, Lemma 2):
 $${\cal J}=\tau^p\mathbb Z\bigoplus_{a=1}^{(p-1)/2}\tau_{-a}\tau^a \mathbb Z.$$
 Observe also that, for $2\leq n\leq p-2,$ we have:
 $$j_{1,p-n}=j_{1,n-1}.$$ Let $V$ be the sub-$\mathbb Z$-module of ${\cal J}$ generated by ${\cal N}$ and the $j_{1,n},$ $1\leq n\leq (p-1)/2.$ Then for $1\leq n\leq p-2,$ $j_{1,n}\in V.$ Furthermore:
 $$\prod_{n=1}^{p-2}j_{1,n}=\frac{\tau^p}{{\cal N}}.$$
 Therefore $\tau ^p \in V.$ Since $\tau_{-1}\tau ^{1}= {\cal N},$ $\tau_{-1}\tau^{1}\in V.$ Now, let $2\leq r \leq (p-1)/2$ and assume that we have proved that $\tau_{-(r-1)}\tau ^{r-1} \in V.$ We have:
 $$j_{1, r-1}=\frac{\tau \tau_{r-1}}{\tau_r}=\frac{ {\cal N} \tau \tau_{1-r}^{-1}}{{\cal N} \tau_{-r}^{-1}}.$$Thus:
 $$\tau_{-r}= j_{1,r-1}^{-1} \tau_{1-r} \tau ^{-1},$$
 and
 $$\tau_{-r}\tau^r = j_{1, r-1}^{-1} \tau_{-(r-1)}\tau^{r-1}.$$
 Thus $\tau_{-r}\tau ^r \in V$ and the Lemma follows. $\diamondsuit$\par
 \newtheorem{Lemma12}[Lemma1]{Lemma}
 \begin{Lemma12}\label{Lemma12}
 Let $\ell$ be a prime number, $\ell\not =p.$ Let $P$ be a prime ideal of $\cal O$ above $\ell$ and let $a\in \{ 1,\cdots , p-2\}.$ Then $\mathbb Q(j_{1,a}(P))=L$ if and only if $\ell \equiv 1\pmod{p}$ and $a^2+a+1\not \equiv 0\pmod{p}$ if $p\equiv 1 \pmod{3}.$
 \end{Lemma12}
 \noindent{\sl Proof}   Since $j_{1,a}(P)\equiv 1\pmod{\pi^2}$ and $j_{1,a}(P)j_{1,a}(P)^{\sigma_{-1}}=\ell^f$ where $f$ is the order of $\ell$ in $(\mathbb Z/p\mathbb Z)^*,$ we have:
 $$\forall \sigma \in \Delta,\,  j_{1,a}(P)^{\sigma }=j_{1,a}(P)\Leftrightarrow j_{1,a}(P)^{\sigma }{\cal{O}}=j_{1,a}(P)\cal{O}.$$
 Recall that:
 $$\forall \sigma \in \Delta,\, j_{1,a}(P)^{\sigma }{\cal{O}}=j_{1,a}(P){\cal{O}}\Leftrightarrow P^{(\sigma -1)(1+\sigma_a-\sigma_{1+a})\theta}=\cal O.$$
Since $j_{1,a}(P)^{\sigma_{\ell }}=j_{1,a}(P),$ we can assume $\ell\equiv 1 \pmod{p}.$ Let $\sigma \in \Delta,$ we have to consider the following equation in $\mathbb C[\Delta ]:$
$$(\sigma-1)(1+\sigma_a-\sigma_{1+a})\theta =0.$$
This is equivalent to:
$$\forall \psi \in {\widehat {\Delta}}, \psi\, {\rm odd},\, (\psi (\sigma )-1)(1+\psi(a)-\psi (1+a))=0.$$
Assume that $\omega^3(\sigma )\not =1.$ Then:
$$1+\omega^3(a)-\omega^3(1+a)=0.$$
This implies:
$$a^2+a\equiv 0\pmod{p},$$
which is a contradiction. Thus $\omega^3(\sigma )=1.$ Let's suppose that $\sigma \not = 1.$  We get:
$$1+\omega (a)=\omega (1+a),$$
which is equivalent to:
$$a^2+a+1\equiv 0\pmod{p}.$$
Conversely, one can see that if $p\equiv 1\pmod{3},$ $a^2+a+1\equiv 0\pmod{p},$ $\omega^3(\sigma )=1,$ then:
$$\forall \psi \in {\widehat {\Delta}}, \psi\, {\rm odd},\, (\psi (\sigma )-1)(1+\psi(a)-\psi (1+a))=0.$$
The Lemma follows. $\diamondsuit$\par

 For $x\in \mathbb Z_p,$ let $[x] \in \{ 0, \cdots , p-1\} $ such that $x\equiv [x] \pmod{p}.$ We set:
 $$\eta=(\prod _{n=1}^{p-2} j_{1,n}^{[n^{-1}]})^{1-\sigma_{-1}} \in{\cal  J}^{-}.$$
 \newtheorem{Lemma3}[Lemma1]{Lemma}
\begin{Lemma3}\label{Lemma3}${}$\par
\noindent a) Let $\psi \in {\widehat \Delta},$ $\psi\not = \omega,$ $\psi$ odd. Then:
$$e_{\psi}(\sum_{n=1}^{p-2} (1+\sigma _n-\sigma_{1+n})[n^{-1}])\in \mathbb Z_p^* e_{\psi}.$$
b) We have:
$$\frac{1}{p}e_{\omega }(\sum_{n=1}^{p-2} (1+\sigma _n-\sigma_{1+n})[n^{-1}])\in \mathbb Z_p^ *  e_{\omega }.$$
\end{Lemma3}
\noindent {\sl Proof} ${}$\par
\noindent a) Write $\psi =\omega ^k,$ $k$ odd, $k\in \{ 3, \cdots , p-2\}.$ We have:
$$\sum_{n=2}^{p-2} (1+\psi(n)-\psi (1+n))[n^{-1}]\equiv \sum_{n=1}^{p-1}\frac{1+n^k-(1+n)^k}{n} \equiv k\pmod{p}.$$
This implies a).\par
\noindent b) We have:
$$\forall a\in \mathbb F_p^*,\, \omega (a) \equiv a^p\pmod{p^2}.$$
Thus:
$$\frac{1}{p}\sum_{n=1}^{p-2}(1+\omega(n)-\omega(1+n))[n^{-1}]\equiv  - \sum_{n=1}^{p-1}\sum_{k=1}^{p-1}\frac {p!}{(p-k)!\, k!\, p }n^{k-1}\pmod{p}.$$
And we get:
$$\frac{1}{p}\sum_{n=1}^{p-2}(1+\omega(n)-\omega(1+n))[n^{-1}]\equiv-1\pmod{p}.$$
This implies b). $\diamondsuit$\par
\newtheorem{Lemma4}[Lemma1]{Lemma}
\begin{Lemma4}\label{Lemma4}
Let $\ell$ be a prime number, $\ell \not = p.$ Let $V_{\ell }$ be the sub-$\mathbb Z[\Delta]$-module of $L^*/(L^*)^p$ generated by $\{ f(P), f\in{\cal J}\}$ where $P$ is some prime of ${\cal I}$ above $\ell .$ Let $\psi \in {\widehat \Delta },$ $\psi$ odd and $\psi \not = \omega.$ Then:
$$V_{\ell }(\psi)= \mathbb F_p e_{\psi} \eta (P).$$
\end{Lemma4}
\noindent {\sl Proof}  Let $E=L(\zeta_{\ell }).$ Then:
$$\frac{L^*}{(L^*)^p}(\psi )\hookrightarrow \frac{E^*}{(E^*)^p}(\psi).$$
Now, in $\frac{E^*}{(E^*)^p}(\psi),$ we have:
$$V_{\ell }(\psi)= \mathbb F_p e_{\psi} \tau(P).$$
It remains to apply Lemma \ref{Lemma3}. $\diamondsuit$\par
Finally, we mention the following Lemma:\par
\newtheorem{Lemma5}[Lemma1]{Lemma}
\begin{Lemma5}\label{Lemma5}
We have:
$$({\cal J}^{-}:\mathbb Z[\Delta ]\eta )=2^{\frac{p-3}{2}}\frac{1}{p} \prod_{\psi \in {\widehat \Delta },\,  \psi \, {\rm odd}}(\sum_{n=1}^{p-2} (1+\psi (n)-\psi (1+n))[n^{-1}]).$$
Furthermore:
$$({\cal J}^{-}:\mathbb Z[\Delta ]\eta )\not \equiv 0\pmod{p}.$$
\end{Lemma5}
\noindent {\sl Proof} Set ${\widetilde {\cal J}^{-}}=(1-\sigma_{-1}){\cal J}\subset {\cal J}^{-}.$ Ten (see \cite{WAS}, paragraph 6.4):
$$({\cal J}^{-}:{\widetilde {\cal J}^{-}})=2^{\frac{p-3}{2}}.$$
Now, by the same kind of arguments as in \cite{WAS}, paragraph 6.4, we get:
$$({\widetilde  {\cal J}^{-}}:\mathbb Z[\Delta ]\eta )=\frac{1}{p} \prod_{\psi \in {\widehat \Delta },\,  \psi \, {\rm odd}}(\sum_{n=1}^{p-2} (1+\psi (n)-\psi (1+n))[n^{-1}]).$$
It remains to apply Lemma \ref{Lemma3}. to conclude the proof of this Lemma. $\diamondsuit$${}$\par

\section {Jacobi Sums and the Ideal Class Group of $\mathbb Q(\zeta_p )$}\par
Recall that the Iwasawa-Leopoldt Conjecture asserts that $A$ is a cyclic $\mathbb Z_p[\Delta ]$-module. This Conjecture is equivalent to:
$$\forall \psi \in {\widehat \Delta}, \psi \, {\rm odd},\psi \not = \omega,\, A(\psi )\simeq \frac{\mathbb Z_p}{B_{1,\psi^{-1}}\mathbb Z_p}.$$
It is well-known (see \cite{WAS}, Theorem 10.9) that:
$$\forall \psi \in {\widehat \Delta}, \psi \, {\rm odd},\psi \not = \omega,\, A(\omega \psi^{-1})=\{ 0\} \Rightarrow A(\psi )\simeq \frac{\mathbb Z_p}{B_{1,\psi^{-1}}\mathbb Z_p}.$$
In this paragraph, we will study the links between Jacobi sums and the structure of $A^{-}.$\par
We fix $\psi \in {\widehat \Delta },$ $\psi$ odd and $\psi \not = \omega.$ We set:
$$m(\psi)=v_p(B_{1,\psi^{-1}}).$$
Recall that, by \cite{WAS}, paragraph 13.6, we have:
$$\mid A(\psi )\mid =p^{m(\psi )}.$$
Let $p^{k(\psi)}$ be the exponent of the group $A(\psi ).$ Then:
$$B_{1,\psi^{-1}}\equiv 0\pmod{p^{k(\psi )}}.$$
\newtheorem{Lemma6}{Lemma}[section]
 \begin{Lemma6} \label{Lemma6}
Let $P$ be a prime ideal in $\cal I$ above a prime number $\ell .$ Then:
$$e_{\psi}\eta (P){\cal O}=0 \, \, {\rm in}\, \frac{{\cal I}}{{\cal I}^p}\, \Leftrightarrow \, \psi(\ell )\not =1\, {\rm or}\, B_{1,\psi^{-1}}\equiv 0\pmod{p}.$$
\end{Lemma6}
\noindent{\sl Proof}
 First note that,  if $\rho \in {\widehat \Delta },$ then $e_{\rho }P=0$  in  $\frac{{\cal I}}{{\cal I}^p}$ if and only if $\rho (\ell )\not =1.$ By the Stickelberger Theorem, we have:
$$\eta (P){\cal O}= (\sum_{n=1}^{p-2} (1+\sigma_n-\sigma_{1+n})[n^{-1}])(1-\sigma_{-1})\theta \, P.$$
Recall that:
$$e_{\psi }\theta =B_{1,\psi^{-1}}e_{\psi }.$$
The Lemma follows. $\diamondsuit$\par
\newtheorem{Lemma7}[Lemma6]{Lemma}
\begin{Lemma7}\label{Lemma7}
Let $f\in {\cal W}^{-}.$ Then $f$ lies in ${\cal W}^p$ if and only if for all prime ideal $P\in {\cal I},$ $f(P)\in (L^*)^p.$
\end{Lemma7}
\noindent{\sl Proof} Let $f\in {\cal W}^{-}$ such that for all prime ideal $P\in {\cal I},$ $f(P)\in (L^*)^p.$ Let ${\cal A}\in {\cal I}.$ Then there exists $\gamma_a\in L^*\cap U$ such that $\gamma_a {\overline \gamma_a}=1$ and:
$$f({\cal A})=\gamma_a^p.$$
Observe that $\beta (f)\in p(\mathbb Z[\Delta ])^{-}.$ Let $g:{\cal I}\rightarrow L^*,$ ${\cal A}\mapsto \gamma_a.$ Then one can verify that $f=g^p$ and $g\in {\cal W}^{-}.$ $\diamondsuit$\par
Let $m\geq 1$ such that $p^m> \mid A\mid .$ Set $n=\mid Cl(L)\mid /\mid A\mid.$ Let $e_m(\psi)\in \mathbb Z[\Delta ]^{-}$ such that:
$$e_m(\psi)\equiv e_{\psi}\pmod{p^m}.$$
Set:
$$\beta_{\psi}=2np^{k(\psi)} e_m(\psi)\in \mathbb Z[\Delta ]^{-}.$$
Since $np^{k(\psi)} e_m(\psi)\in ({\rm Ann}_{\mathbb Z[\Delta }Cl(L))^{-},$ by Lemma \ref{Lemma1}, there exists a unique element $f_{\psi}\in {\cal W}^{-}$ such that $\beta (f_{\psi})=\beta_{\psi }.$ Recall that:
$$({\rm Ann}_{\mathbb Z_p[\Delta ]} A)(\psi )=p^{k(\psi )} \mathbb Z_p e_{\psi}.$$
Therefore, for $0\leq k\leq m,$  $\frac{{\cal W}^{-}}{({\cal W}^{-})^{p^k}}(\psi )$ is cyclic of order $p^k$ generated by the image of $f_{\psi}.$ We set:
$$W=\{ f({\cal A}),\, {\cal A}\in {\cal I}, f\in {\cal W}\},$$
and:
$$J=\{ f({\cal A}),\, {\cal A}\in {\cal I}, f\in {\cal J}\}.$$
Observe that $J$ is a sub-$\mathbb Z[\Delta ]$-module of $W,$ and it is called the module of Jacobi sums of $\mathbb Q(\zeta_p).$ Note that, by Lemma \ref{Lemma7}, we have:
$$\frac{W(L^{*})^p}{(L^*)^p}(\psi )\not = \{ 0 \} .$$
\newtheorem{Theorem1}[Lemma6]{Theorem}
\begin{Theorem1}\label{Theorem1}
The map $f_{\psi}$ induces an isomorphism of groups:
$$A(\psi)\simeq \frac{W(L^*)^{p^{k(\psi)}}}{(L^*)^{p^{k(\psi)}}}(\psi ).$$
\end{Theorem1}
\noindent{\sl Proof}  First observe that $m\geq k(\psi )+1.$ Let $P$ be a prime in $\cal I.$ Then:
$$f_{\psi}(P){\cal O}=P^{\beta_{\psi}}.$$
Let $\rho \in {\widehat \Delta}, $ $\rho \not = \psi.$ Then:
$$e_m(\rho)e_m(\psi) \equiv 0\pmod{p^m}.$$
Therefore, there exists $\gamma \in L^*\cap U$ such that:
$$P^{(1-\sigma_{-1})ne_m(\rho)e_m(\psi)}=(\frac {\gamma}{\sigma_{-1} (\gamma)})^p{\cal O}.$$
But $(1-\sigma_{-1})e_m(\psi)=2e_m(\psi).$ Thus, there exists $\alpha \in L^*\cap U,$ $\alpha \sigma_{-1}(\alpha)=1,$ and:
$$f_{\psi}(P)^{e_m(\rho)}=\alpha ^{p^{k(\psi)+1}}.$$
Therefore, $e_{\rho} f_{\psi}({\cal I})=0$ in $L^*/(L^*)^{p^{k(\psi )+1}}.$\par
\noindent It is clear that $f_{\psi}$ induces a morphism:
$$\frac{{\cal I}}{({\cal I})^{p^m}{\cal P}}(\psi )\rightarrow \frac {L^*}{(L^*)^{p^{k(\psi)}}}(\psi ).$$
Now, let $P$ be a prime in $\cal I$ such that $e_{\psi}f_{\psi}(P)=0$ in $\frac {L^*}{(L^*)^{p^{k(\psi)}}}(\psi ).$ Then, by the above remark, we get:
$$f_{\psi}(P)=0\, {\rm in}\,  L^*/(L^*)^{p^{k(\psi)}}.$$
Thus, there exists $\gamma \in L^*\cap U$ such that:
$$P^{\beta_{\psi}}=\gamma ^{p^{k(\psi)}}{\cal O}.$$
Thus:
$$P^{2ne_m(\psi)}=\gamma {\cal O}.$$
This implies:
$$e_{\psi} P=0\, {\rm in}\, \frac{{\cal I}}{({\cal I})^{p^m}{\cal P}}(\psi ).$$
Thus our map is injective. Now, observe that the image of the map induced by $f_{\psi}$ is  $\frac{W(L^*)^{p^{k(\psi)}}}{(L^*)^{p^{k(\psi)}}}(\psi )$ and that:
$$A(\psi )\simeq \frac{{\cal I}}{({\cal I})^{p^m}{\cal P}}(\psi ).$$
The Theorem follows. $\diamondsuit$\par
Recall that:
$$\eta = (\prod_{n=1}^{p-2} j_{1,n}^{[n^{-1}]})^{1-\sigma_{-1}}\in {\cal J}^{-}.$$
Set:
$$z=(1-\sigma_{-1})\sum_{n=1}^{p-2}(1+\sigma_n -\sigma_{1+n})[n^{-1}]\, \in \mathbb Z[\Delta ]^{-}.$$
We have:
$$\beta(\eta )=z\theta .$$
\newtheorem{Corollary1}[Lemma6]{Corollary}
\begin{Corollary1} \label{Corollary1}${}$\par
\noindent 1) The map $\eta $ induces an isomorphism of groups:
$$A(\psi)\simeq \frac{J(L^*)^{p^{m(\psi )}}}{(L^*)^{p^{m(\psi )}}}(\psi ).$$
2) $\frac {J (L^* )^p}{(L^*)^p}(\psi )\not = \{ 0\}$ if and only if $A(\psi )$ is $\mathbb Z_p$-cyclic.
\end{Corollary1}
\noindent {\sl Proof} ${}$\par
\noindent 1) Let $P$ be a prime in $\cal I.$ Then one can show that:
$$f_{\psi }(P)^{z\theta }=\eta (P)^{2np^{k(\psi )}e_m(\psi )}.$$
The first assertion follows from Theorem \ref{Theorem1}.\par
\noindent 2) Note that:
$$A(\psi)\, {\rm is}\, \mathbb Z_p{\rm -cyclic}\, \Leftrightarrow m(\psi)=k(\psi ).$$
Thus, if $A(\psi)$ is $\mathbb Z_p$-cyclic, then:
$$\frac {J (L^* )^p}{(L^*)^p}(\psi )=\frac {W (L^* )^p}{(L^*)^p}(\psi )\not = \{ 0\}.$$
By the proof of assertion 1), if $k(\psi )<m(\psi)$ and if $P$ is a prime in $\cal I,$ then:
$$\eta(P)^{e_m(\psi )}\in (L^*)^p.$$
Therefore, we get assertion 2). $\diamondsuit$\par


\section{The $p$-adic behavior of Jacobi Sums}\par
Let $M$ be a subgroup of $L^*/(L^*)^p,$ we  say that  $M$ is unramified if $L({}^p\sqrt{M})/L$ is an unramified extension. Note that Kummer's Lemma asserts that (\cite{WAS}, Theorem 5.36):
$$\forall \rho \in {\widehat \Delta},\, \rho \, {\rm even},\, \rho \not =1,\, 
\frac{E}{E^p}(\rho)\, {\rm is }\, {\rm unramified } \Rightarrow B_{1,\rho \omega^{-1}}\equiv 0\pmod{p}.$$
It is natural to ask if  this implication is in fact an equivalence (see \cite{ANG}, \cite{ASS&NGU}). We will say that the converse of Kummer's Lemma is true for the character $\rho $ if we have:
$$\frac{E}{E^p}(\rho)\, {\rm is }\, {\rm unramified } \Leftrightarrow B_{1,\rho \omega^{-1}}\equiv 0\pmod{p}.$$
 In this paragraph, we will study this question with the help of Jacobi sums.\par
 Let $F/L$ be the maximal abelian $p$-extension of $L$ which is unramified outside $p.$ Set ${\cal X}={\rm Gal}(F/L).$ We have an exact sequence of $\mathbb Z_p[\Delta]$-modules (\cite{WAS}, Corollary 13.6):
 $$0\rightarrow \frac{U}{\overline E}\rightarrow {\cal X}\rightarrow A\rightarrow 0\, .$$
Let $\rho \in {\widehat \Delta},$ observe that:\par
\noindent - if $\rho =1, \omega $ then ${\cal X}(\rho )\simeq \mathbb Z_p,$\par
\noindent - if $\rho$ is even,$ \rho \not = 1,$ ${\cal X}(\rho)\simeq {\rm Tor}_{\mathbb Z_p}{\cal X}(\rho),$\par
\noindent - if $\rho$ is odd, $\rho \not = \omega,$ ${\cal X}(\rho)\simeq \mathbb Z_p\oplus {\rm Tor}_{\mathbb Z_p}{\cal X}(\rho).$\par
\newtheorem{Lemma8}{Lemma}[section]
\begin{Lemma8}\label{Lemma8}
Let $\psi \in {\widehat \Delta},$ $\psi$ odd, $\psi \not =\omega .$ Then:
$$d_p{\rm Tor}_{\mathbb Z_p}{\cal X}(\psi )=d_pA(\omega \psi^{-1}).$$
\end{Lemma8}
\noindent{\sl Proof}
This is a consequence of the proof of Leopoldt's reflection Theorem (\cite{WAS}, Theorem 10.9). For the convenience of the reader, we give the proof of the above Lemma.\par
Let $H$ be the Galois group of the maximal abelian extension of $L$ which is unramified outside $p$ and of exponent $p.$ Then $H$ is a $\mathbb Z_p[\Delta ]$-module and we have:\par
\noindent - $H(1)\simeq \mathbb F_p$ and corresponds to $L(\zeta_{p^2})/L,$\par
\noindent - $H(\omega )\simeq \mathbb F_p$ and corresponds to $L({}^p\sqrt{p})/L,$\par
\noindent - if $\rho$ is even, $\rho \not =1,$ $d_pH(\rho)=d_p{\rm Tor}_{\mathbb Z_p}{\cal X}(\rho),$\par
\noindent- if $\rho$ is odd,$ \rho \not = \omega,$ then $d_pH(\rho)=1+d_p{\rm Tor}_{\mathbb Z_p}{\cal X}(\rho).$\par
\noindent Let $V$ be the sub-$\mathbb Z[\Delta ]$-module of $L^*/(L^*)^p$ which corresponds to $H,$ i.e. $H={\rm Gal}(L({}^p\sqrt{V})/L).$ Let $M$ be the sub-$\mathbb Z[\Delta ]$-module of $L^*/(L^*)^p$ generated by $E$ and $1-\zeta_p .$ We have an exact sequence:
$$0\rightarrow M\rightarrow V\rightarrow A[p]\rightarrow 0.$$
Soit $\psi \in \widehat \Delta,$ $\psi $ odd, $\psi \not = \omega .$ We have by Kummer theory :
$$1+d_p{\rm Tor}_{\mathbb Z_p}{\cal X}(\psi)=d_pV(\omega \psi^{-1}).$$
And, by the above exact sequence, we have:
$$d_pV(\omega \psi^{-1})=1+d_pA(\omega\psi^{-1}).$$
The Lemma follows. $\diamondsuit$\par
\newtheorem{Lemma9}[Lemma8]{Lemma}
\begin{Lemma9}\label{Lemma9}
Let $\rho \in \widehat \Delta,$ $\rho$ even and $\rho \not =1.$ If  $\frac{E}{E^p}(\rho)$ is ramified then $d_pA(\rho )=d_pA(\omega \rho^{-1}).$
\end{Lemma9}
\noindent{\sl Proof}  We keep the notations  of the proof of Lemma \ref{Lemma8}. Let $V^{nr}\subset V$ which corresponds via Kummmer theory to $A/pA.$ Then:
$$V^{nr}(\rho) \simeq \frac{A}{pA}(\omega \rho^{-1}).$$
But we have:\par 
\noindent $\frac{E}{E^p}(\rho)$ is ramified if and only if $V^{nr}(\rho)\hookrightarrow A[p](\rho ).$\par
\noindent Now recall that:
$$d_p A(\rho )\leq d_p A(\omega \rho^{-1}).$$
The Lemma follows. $\diamondsuit$\par
\newtheorem{Lemma10}[Lemma8]{Lemma}
\begin{Lemma10}\label{Lemma10}
There exists an unique $\mathbb Z[\Delta ]$-morphism $\varphi: K^*\rightarrow \mathbb Z_p[\Delta ]$ such that:
$$\forall x\in K^*,\, \varphi (x)\zeta_p = {\rm Log}_p(x).$$
Furthermore, we have:
$${\rm Im}\, \varphi=\bigoplus_{\rho =1,\omega}p\mathbb Z_pe_{\rho}\, \bigoplus_{\rho\not = 1,\omega}\mathbb Z_p e_{\rho }.$$
\end{Lemma10}
\noindent{\sl Proof} Let $\lambda \in K^*$ such that $\lambda^{p-1}=-p.$ Then:
$$K^*=\lambda^{\mathbb Z}\times \mu_{p-1}\times \mu_p \times U.$$Recall that:\par
\noindent - the kernel of ${\rm Log}_p$ on $K^*$ is equal to $\lambda^{\mathbb Z}\times \mu_{p-1}\times \mu_p,$\par
\noindent - ${\rm Log}_p U= \pi^2\mathbb Z_p[\zeta_p].$\par
\noindent For $\rho \in \widehat \Delta,$ set:
$$\tau (\rho ) = \sum_{a=1}^{p-1} \rho (a) \zeta_p‰ \in \mathbb Z_p[\zeta_p].$$
Then:
$$e_{\rho }\zeta_p =\tau (\rho^{-1}).$$
But recall that that $\mathbb Z_p[\zeta_p]=\mathbb Z_p[\Delta ] \zeta_p.$ Thus:
$$e_{\rho }\mathbb Z_p[\zeta_p]=\mathbb Z_p \tau (\rho^{-1}).$$
If $\rho =\omega ^k,$ $k\in \{ 0\cdots, p-2\},$ we have:
$$v_p(\tau (\rho^{-1}))=\frac{k}{p-1}.$$
Therefore:
$$\pi^2\mathbb Z_p[\zeta_p]=\bigoplus_{\rho=1,\omega } p\mathbb Z_p \tau(\rho^{-1}) \, \bigoplus _{\rho \not =1, \omega }\mathbb Z_p \tau (\rho^{-1}).$$
The Lemma follows. $\diamondsuit$\par
Let $P$ be a prime in $\cal I.$ We fix a generator $r_P\in \mathbb F_P^*$ such that:
$$\chi_P(r_P)=\zeta_p.$$
For $x\in \mathbb F_P^*,$ let ${\rm Ind}(P,x)\in \{ 0,\cdots ,NP-2\}$ such that:
$$x=r_p^{{\rm Ind}(P,x)}.$$
We recall the following Theorem (see also \cite{THA} for a statement similar but weaker than   part 2) of the following Theorem):
\newtheorem{Theorem2}[Lemma8]{Theorem}
\begin{Theorem2}\label{Theorem2}
${}$\par
\noindent 1) $\varphi (1-\zeta_p)=\sum_{\rho \in {\widehat \Delta}, \rho \not =1, \rho\, {\rm even}} -(p-1)^{-1} L_p(1,\rho ) e_{\rho }.$\par
\noindent 2) Let $\psi \in \widehat \Delta,$ $\psi $ odd, $\psi \not = \omega.$ Write $\psi =\omega^k,$ $k\in \{ 2,\cdots, p-2\}.$ We have:
$$e_{\psi }\varphi (\eta (P))\equiv 2k\, {\rm Ind}(P, \prod_{a=1}^{p-1}(\frac{1-\zeta_p^{-a}}{1-\zeta_p})^{a^{k-1}})\, e_{\psi}\pmod{p}.$$
\end{Theorem2}
\noindent{\sl Proof}${}$\par
\noindent 1) Let $\rho \in \widehat \Delta,$ $\rho$ even, $\rho \not =1.$ By \cite{WAS}, Theorem 5.18, we have:
$$L_p(1,\rho)\tau(\rho^{-1})=-(p-1) e_{\rho} {\rm Log}_p(1-\zeta_p).$$
Thus the first assertion follows.\par
\noindent 2) Let $\psi \in \widehat \Delta,$ $\psi $ odd, $\psi \not = \omega.$
By a beautiful result of Uehara (\cite{UEH}, Theorem1), we have:
$$e_{\psi}{\rm Log}_p (\eta (P))\equiv 2k\, {\rm Ind}(P, \prod_{a=1}^{p-1}(\frac{1-\zeta_p^{-a}}{1-\zeta_p})^{a^{k-1}})\, \tau (\psi^{-1}) \pmod{p}.$$
This implies the second assertion. $\diamondsuit$ \par
\newtheorem{Theorem3}[Lemma8]{Theorem}
\begin{Theorem3}\label{Theorem3}
Let $\psi \in \widehat \Delta,$ $\psi \not = \omega,$ $\psi$ odd. We have the following exact sequences:
$$0\rightarrow {\rm Tor}_{\mathbb Z_p} {\cal X}(\psi)\rightarrow A(\psi )\rightarrow \frac{{\overline W}(\psi )}{U^{p^{k(\psi)}}(\psi )}\rightarrow 0,$$
$$0\rightarrow {\rm Tor}_{\mathbb Z_p} {\cal X}(\psi)\rightarrow A(\psi )\rightarrow \frac{{\overline J}(\psi )}{U^{p^{m(\psi)}}(\psi )}\rightarrow 0.$$
\end{Theorem3}
\noindent{\sl Proof} This Theorem is a consequence of the method developped by Iwasawa in \cite{IWA1}. Let's recall briefly this method.\par
Let $f\in {\cal W}.$ Set, for $n\geq 2,$ ${\cal P}_n=\{ \alpha {\cal O},\, \alpha \equiv 1\pmod{\pi^n}\}.$ Observe that:
$$f({\cal P}_n)\subset 1+\pi^n \mathbb Z_p[\zeta_p].$$
Let:
$${\widetilde {\cal X}}=\lim_{\leftarrow}\frac{{\cal I}}{{\cal P}_n}.$$
If $\widetilde F$ is the maximal abelian extension of $L$ which is unramified outside $p,$ then, we get by class field theory:
$${\widetilde {\cal X}}\simeq {\rm Gal}({\widetilde F}/L).$$
By \cite{WAS}, Theorem 13.4, the natural surjective map ${\widetilde {\cal X}}\rightarrow {\cal X}$ has a finite kernel of order prime to $p.$ Thus $f$ induces a map:
$${\overline f}: {\cal X}\rightarrow U.$$
Furthermore:
$${\overline f}(U) =U^{\beta (f)}\subset {\overline f}({\cal X}).$$
Now let $\psi \in  \widehat \Delta,$ $\psi $ odd, $\psi \not = \omega.$ We have a map:
$${\overline f}: {\cal X}(\psi) \rightarrow U(\psi ).$$
But:
$${\cal X}(\psi) \simeq \mathbb Z_p\bigoplus {\rm Tor}_{\mathbb Z_p} {\cal X}(\psi ),$$
and:
$$U(\psi )\simeq \mathbb Z_p.$$
Thus, if $e_{\psi }\beta (f) \not =0,$ we get:
$${\rm Ker}\, ({\overline f}: {\cal X}(\psi) \rightarrow U(\psi ))={\rm Tor}_{\mathbb Z_p} {\cal X}(\psi ).$$
Therefore, if $e_{\psi }\beta (f) \not =0,$ we get the following exact sequence induced by $f$:
$$0\rightarrow {\rm Tor}_{\mathbb Z_p} {\cal X}(\psi )\rightarrow A(\psi )\rightarrow \frac{{\overline f}({\cal X})(\psi )}{U^{\beta (f)}(\psi )}\rightarrow 0.$$
It remains to apply this construction to $f_{\psi}$ and $\eta $ to get the desired exact sequences. $\diamondsuit$\par
\newtheorem{Corollary2}[Lemma8]{Corollary}
\begin{Corollary2} \label{Corollary2}
${}$\par
\noindent 1) Let $\psi \in \widehat \Delta,$ $\psi $ odd, $\psi \not = \omega.$ Then:
$$d_pA(\psi )=1+d_pA(\omega \psi ^{-1})\Leftrightarrow B_{1,\psi ^{-1}}\equiv 0\pmod{p}\, {\rm and}\, {\overline W}(\psi)=U(\psi ).$$
2) Let $\rho \in \widehat \Delta,$ $\rho $ even and $\rho \not =1.$ Assume that $B_{1,\rho \omega^{-1}}\equiv 0\pmod{p}$ and that ${\overline W}(\omega \rho^{-1})=U(\omega \rho^{-1})$ then the converse of Kummer's Lemma is true for the character $\rho.$
\end{Corollary2}
\noindent {\sl Proof}${}$\par
\noindent 1) We apply Theorem \ref{Theorem3}. We identify  ${\rm Tor}_{\mathbb Z_p} {\cal X}(\psi )$ with its image in $A(\psi ).$ We can write:
$$A(\psi )=B\bigoplus C,$$
where $C$ is cyclic of order $p^{k(\psi)}$ and $B\subset {\rm Tor}_{\mathbb Z_p} {\cal X}(\psi ).$
Now:
$$(C:C\cap {\rm Tor}_{\mathbb Z_p} {\cal X}(\psi ))=({\overline W}(\psi ): U^{p^{k(\psi )}}(\psi )).$$
It remains to apply Lemma \ref{Lemma8} to get the desired result.\par
\noindent 2) We apply the first assertion and Lemma \ref{Lemma8}, we get:
$$d_p A(\rho )=d_pA (\omega \rho^{-1}) -1.$$
It remains to apply Lemma \ref{Lemma9}. $\diamondsuit$\par
We set:
$$W^{nr}=\{ \alpha \in W,\, \alpha \in U^p\}.$$
Let $\psi \in \widehat \Delta,$ $\psi $ odd, $\psi \not =\omega.$ We assume that $B_{1,\psi^{-1}}\equiv 0\pmod{p}.$ Write:
$$A(\psi )=\frac{\mathbb Z}{ p^{e_1}\mathbb Z}\bigoplus \cdots \bigoplus \frac{\mathbb Z}{ p^{e_t}\mathbb Z},$$
where  $t=d_pA(\psi )$ and $1\leq e_1\leq \cdots \leq e_t=k(\psi ).$ 
Set:
$$n(\psi )=\mid\{ i\in \{ 1,\cdots t\},\, e_i=k(\psi )\} \mid.$$
\newtheorem{Corollary3}[Lemma8]{Corollary}
\begin{Corollary3}\label{Corollary3}
We have:
$$n(\psi )-1\leq {\rm dim}_{\mathbb F_p}\frac{W^{nr}(L^*)^p}{(L^*)^p}\leq n(\psi ).$$
Furthermore:
$${\rm dim}_{\mathbb F_p}\frac{W^{nr}(L^*)^p}{(L^*)^p}=n(\psi )\Leftrightarrow {\overline W}(\psi )\not = U(\psi ).$$
\end{Corollary3}
\noindent{\sl Proof} By Theorem\ref{Theorem3} and Theorem \ref{Theorem1}, we have:
$$\frac{W^{nr}(L^*)^{p^{k(\psi )}}}{(L^*)^{p^{k(\psi )}}}\simeq{\rm Ker}(A(\psi )\rightarrow  \frac{{\overline W}(\psi )}{U^{p^{k(\psi)}}(\psi )}).$$
The Corollary follows. $\diamondsuit$\par
\newtheorem{Corollary4}[Lemma8]{Corollary}
\begin{Corollary4}\label{Corollary4}
Assume that $pA^{-}=\{0\}.$ Then we have an isomorphism of groups:
$${\rm Gal}(L({}^p\sqrt{W^{nr}})/L)\simeq \frac{A^{+}}{pA^{+}}.$$
\end{Corollary4}
\noindent{\sl Proof} This result is a consequence of Kummer theory, Corollary \ref{Corollary3} and Corollary \ref{Corollary2}. $\diamondsuit$\par
Note that the above results lead us to ask the following problem (which is a restatement of the converse of Kummer's Lemma):\par
\centerline{do we have $\varphi ({\overline W}^{-})=({\rm Im}\, \varphi )^{-}?$}\par
\noindent Observe  that $e_{\omega}\varphi ({\overline W}^{-})=e_{\omega}({\rm Im}\, \varphi )^{-},$
and since $K_4(\mathbb Z)=\{ 0\},$ we have $A(\omega^{-2})=\{0\}$ (see \cite{KUR}) and therefore $
e_{\omega^3}\varphi ({\overline W}^{-})=e_{\omega^3}({\rm Im}\, \varphi )^{-}.$\par

\section{Remarks on the Jacobian of the Fermat Curve over a finite field}\par
First we fix  some notations and recall some basic facts about global function fields. \par
Let $\mathbb F_q$ be a finite field having $q$ elements. Let $\ell$ be the charactersitic of $\mathbb F_q,$ $\ell \not =p.$ Let $\overline{\mathbb F_q}$ be a fixed algebraic closure of $\mathbb F_q$ and let $\widetilde{\mathbb F_q}= \cup _{n\geq 1,\, n\not \equiv 0\pmod{p}}\mathbb F_{q^n}\subset \overline{\mathbb F_q}.$ Let $k/\mathbb F_q$ be a global function field such that $\mathbb F_q$ is algebraically closed in $k.$ We set:\par
\noindent - $D_k:$ the group of divisors of $k,$\par
\noindent - $D_k^{0}:$ the group of divisors of degree zero of $k,$\par
\noindent - $P_k:$ the group of principal divisors of $k,$\par
\noindent - $J_k:$ the jacobian of $k,$ note that we have:
$$\forall n\geq 1,\, J_k(\mathbb F_{q^n})\simeq \frac{D_{\mathbb F_{q^n}k}^{0}}{P_{\mathbb F_{q^n}k}},$$
- $g_k:$ the genus of $k,$\par
\noindent - $L_k(Z)\in \mathbb Z[Z]:$ the numerator of the zeta function of $k,$ we recall that:
$$\frac{L_k(Z)}{(1-Z)(1-qZ)}=\prod_{v\, {\rm place \, of\, }k}(1-Z^{deg\, v})^{-1},$$
furthermore ${\rm deg}_Z L_k(Z)=2g_k$ and $L_k(1)=\mid J_k(\mathbb F_q)\mid,$\par
\noindent - $C_k(\mathbb F_{q^n})=J_k(\mathbb F_{q^n})\otimes_{\mathbb Z}\mathbb Z_p,$\par
\noindent - $\widetilde{d}_p J_k=d_pC_k(\widetilde{\mathbb F_q}),$ observe that there exists an integer  $m,$ $m\not \equiv 0\pmod{p},$ such that :
$$C_k(\widetilde{\mathbb F_q})=C_k(\mathbb F_{q^m}).$$
Write:  $$L_k(Z)=\prod_{i=1}^{2g_k}(1-\alpha_i Z).$$
For simplicity, we assume that $v_p(\alpha_i -1)>0$ for $i=1,\cdots , 2g_k.$ In this case, we have:
$$C_k(\widetilde{\mathbb F_q})=C_k(\mathbb F_q).$$
Set $P_k(Z)=\prod_{i=1}^{2g_k}(Z-(\alpha_i-1)).$ 
 Let $\gamma$ be the Frobenius of $\mathbb F_q,$ and set:
$$C_n(k)=C_k(\mathbb F_{q^{p^n}}).$$
Let $C_{\infty }(k)=\cup_{n\geq 0} C_n(k),$ and set:
$$M_k={\rm Hom}(\frac{\mathbb Q_p}{\mathbb Z_p}, C_{\infty}(k)).$$
Then $M_k$ is isomorphic to the $p$-adic Tate-module of $J_k.$ Set $\Lambda =\mathbb Z_p[[Z]]$ where $Z$ corresponds to $\gamma -1.$ Then it is well-known that:\par
\noindent - $M_k$ is a $\Lambda$-module of finite type and of torsion,\par
\noindent - as a $\mathbb Z_p$-module $M_k$ is isomorphic to $\mathbb Z_p^{2g_k},$\par
\noindent - $M_k/\omega_nM_k\simeq C_n(k),$ where $\omega_n =(1+Z)^{p^n}-1,$\par
\noindent - ${\rm Char}_{\Lambda}M_k=P_k(Z)\Lambda,$\par
\noindent - the action of $Z$ on $M_k$ is semi-simple, i.e. the minimal polynomial of the action of $Z$ on $M_k$ has only simple roots.\par
Now, let $\ell$ be a prime number, $\ell \not =p.$ We fix a prime $P$ of $\cal O$ above $\ell $ and we view ${\cal O}/P$ as a subfield of $\overline {\mathbb F_{\ell}},$ thus $\mathbb F_q=\cal {O}/P \subset \widetilde{\mathbb F_{\ell}}.$ We identify $\zeta_p$ with its image in $\mathbb F_q.$ Let $X$ be an indeterminate over $\mathbb F_q,$ we set $k=\mathbb F_{\ell}(X,Y)$ where $X^p+Y^p=1,$ and we set: $T=X^p.$ For $a,b\in \mathbb Z,$ let $\tau_{a,b}\in{\rm Gal}(\overline{\mathbb F_{\ell}}k/\overline{\mathbb F_{\ell}}(T))$ such that:
$$\tau_{a,b}(X)=\zeta_p^a  X\, {\rm and}\, \tau_{a,b}(Y)=\zeta_p^bY.$$
Let $a\in\{ 1,\cdots ,p-2\}.$ Let $H_a$ be the subgroup of ${\rm Gal}(\overline{\mathbb F_{\ell}}k/\overline{\mathbb F_{\ell}}(T))$ generated by $\tau_{1, [-a^{-1}]}.$ Set:
$$E_a=\mathbb F_{\ell}(T, XY^{a}).$$
We set:
$$E=\mathbb F_qE_a,$$
$$F=\mathbb F_q k,$$
and observe that:
$$\widetilde{\mathbb F_{\ell}}=\widetilde {\mathbb F_q}.$$
It is clear that:
$$F^{H_a}=E.$$
Finally, we set:
$$G={\rm Gal}(E/\mathbb F_q(T)).$$
Note that $g_E=(p-1)/2.$\par
\newtheorem{Lemma11}{Lemma}[section]
\begin{Lemma11}\label{Lemma11}
We have:
$$L_E(Z)=\prod_{\sigma \in \Delta }(1-j_{1,a}(P)^{\sigma }Z).$$
\end{Lemma11}
\noindent{\sl Proof} Let $\chi \in \widehat G$ such that:
$$\chi (g)=\zeta_p^{-1},$$
where $g\in G$ is such that $g(XY^{a})=\zeta_p XY^{a}.$
Note that:
$$L_E(Z)=\prod_{\sigma \in \Delta }L(Z, \chi^{\sigma }),$$
where:
$$L(Z,\chi )=\prod_{v\,{\rm place\, of\,}\mathbb F_q(T)}(1-\chi (v)Z^{deg\, v})^{-1}.$$
Since $2g_e=p-1,$ we get:
$${\rm deg}_ZL(Z,\chi )=1.$$
For $b\in \mathbb F_q\setminus \{ 0,1\},$ we denote the Frobenius of $T-b$ in $E/\mathbb F_q(T)$ by $Frob_b.$ We have:
$$Frob_b(XY^{a})=(b(1-b)^{a})^{(q-1)/p}XY^{a}.$$
But:
$$L(Z,\chi) \equiv 1+(\sum_{b\in \mathbb F_q\setminus\{ 0,1\}} \chi (Frob_b))X\pmod {X^2}.$$
Thus:
$$L(Z,\chi)=1+(\sum_{b\in \mathbb F_q\setminus\{ 0,1\}} \chi (Frob_b))X.$$
But, we can write:
$$j_{1,a}(P)=-\sum_{i=0}^{p-1} N_i \zeta_p^{-i},$$
where $N_i=\mid \{ \alpha \in \mathbb F_q\setminus \{ 0,1\},\, (\alpha (1-\alpha )^{a})^{(q-1)/p}\equiv \zeta_p^{-i}\pmod{P}\}\mid.$ Therefore:
$$j_{1,a}(P)=-\sum_{b\in \mathbb F_q\setminus \{ 0,1\} }\chi (Frob_b).$$
The Lemma follows. $\diamondsuit$\par
\newtheorem{Theorem4}[Lemma11]{Theorem}
\begin{Theorem4} \label{Theorem4}
Let $n$ be the smallest integer (if it exists) such that $3\leq n\leq p-2,$ $n$ odd and $e_{\omega^n}j_{1,a}(P) \not \in U^p,$ then:
$$J_k(\widetilde{\mathbb F_{\ell}})^{H_a}\otimes_{\mathbb Z} \mathbb Z_p\simeq (\frac{\mathbb Z}{p \mathbb Z})^{n}.$$
If such an integer doesn't exist  then:\par
\noindent 1) $\widetilde{d}_p J_k^{H_a}=p-1,$\par
\noindent 2) we have:
$$J_k(\widetilde{\mathbb F_{\ell}})^{H_a}\otimes_{\mathbb Z} \mathbb Z_p\simeq (\frac{\mathbb Z}{p \mathbb Z})^{p-1}\, \Leftrightarrow \ell ^{p-1}\not \equiv 1\pmod{p^2}.$$
\end{Theorem4}
\noindent{\sl Proof} The proof of this result is based on ideas developped by Greenberg in \cite{GRE}.  Set $H=H_a.$ Let $P_0$ be the prime of $E$ above $T,$ $P_1$ the prime of $E$ above $T-1$ and $P_{\infty}$ the prime of $E$ above $\frac{1}{T}.$ Recall that we have in $D_E:$
$$p(P_0-P_{\infty })=(T),$$
$$p(P_1-P_{\infty})= (T-1),$$
$$P_0- P_{\infty} + a(P_1-P_{\infty })=(XY^{a}).$$
Thus, by \cite{GRE}, paragraph 2, we get:
$$J_E(\mathbb F_q)^G\simeq \frac{\mathbb Z}{p\mathbb Z},$$
and $J_E(\mathbb F_q)^G$ is generated by the class of $P_0-P_{\infty}.$ Observe also that $F/E$ is unramified and cyclic of order $p.$ Let's start by the folowing exact sequence:
$$0\rightarrow \mathbb F_q^*\rightarrow F^*\rightarrow P_F\rightarrow 0.$$
We get:
$$\frac{P_F^{H}}{P_E}\simeq \frac {\mathbb Z}{p\mathbb Z},$$
and $\frac{P_F^{H}}{P_E}$ is generated by the image of $P_0-P_{\infty}$ in $D_F.$ In particular:
$$\frac{P_F^{H}}{P_E}\simeq J_E(\mathbb F_q)^G.$$
Note that we also have:
$$0\rightarrow H^1(H,P_F)\rightarrow H^2(H,\mathbb F_q^*)\rightarrow H^2(H, F^*).$$
But $F/E$ is unramified and cyclic, therefore every element of $\mathbb F_q^*$ is a norm  in the extension $F/E.$ Thus:
$$H^1(H,P_F)\simeq \frac{\mathbb Z}{p\mathbb Z}.$$
Now, we look at the exact sequence:
$$0\rightarrow P_F\rightarrow D_F^0\rightarrow J_F(\mathbb F_q)\rightarrow 0.$$
Since $F/E$ is unramified:
$$H^1(H, D_F^0)=\{ 0\}.$$
Therefore, we have obtained the following exact sequence:
$$0\rightarrow J_E(\mathbb F_q)^G\rightarrow J_E(\mathbb F_q)\rightarrow J_F(\mathbb F_q)^H\rightarrow \frac{\mathbb Z}{p\mathbb Z}\rightarrow 0.$$
Now, it is not difficult to deduce that, for all $n\geq 1,$ we have the following exact sequence:
$$0\rightarrow \frac{\mathbb Z}{p\mathbb Z}\rightarrow J_E(\mathbb F_{q^n})\rightarrow J_F(\mathbb F_{q^n})^H\rightarrow \frac{\mathbb Z}{p\mathbb Z}\rightarrow 0.$$
From this, we get the following exact sequence of $\mathbb Z_p[G]$-modules and $\Lambda$-modules:
$$0\rightarrow M_E\rightarrow M_F^{H}\rightarrow \frac{\mathbb Z}{p\mathbb Z}\rightarrow 0.$$
Recall that in our situation, by Lemma \ref{Lemma11}, we have:
$$P_E(Z)=\prod_{\sigma \in \Delta}(Z-(j_{1,a}(P)^{\sigma }-1)).$$
Furthermore the action of $G$ and $Z$ commute on $M_F^H.$
Now, we have:\par
\noindent - ${\rm Char}_{\Lambda} M_F^H={\rm Char}_{\Lambda}M_E=P_E(Z)\Lambda ,$\par
\noindent - $M_F^H\simeq \mathbb Z_p^{p-1}$ as $\mathbb Z_p$-modules,\par
\noindent - $M_F^H/\omega_n\simeq C_n(F)^H.$\par
\noindent Observe that:
$$C_0(F)^H=J_k(\widetilde{\mathbb F_{\ell}})^{H_a}\otimes_{\mathbb Z}\mathbb Z_p.$$
Note also that the minimal polynomial of the action of $Z$ on $M_F^H$ is:
$${\rm Irr}(j_{1,a}(P)-1, \mathbb Q_p; Z):=G(Z).$$
Set $N=\sum_{\delta \in G}\delta.$ Then one can see that:
$$NM_E=NM_F^H=\{ 0\}.$$ Thus $M_F^H$ is a $\mathbb Z_p[G]/N\mathbb Z_p[G]$-module. Now, we identify $\mathbb Z_p[G]/N\mathbb Z_p[G]$ with $\mathbb Z_p[\zeta_p].$ Since $M_F^H\simeq \mathbb Z_p^{p-1},$ there exists $m\in M_F^H$ such that:
$$M_F^H\simeq \mathbb Z_p[\zeta_p].m,$$
i.e. $M_F^H$ is  a free $\mathbb Z_p[\zeta_p]$-module of rank one. Therefore there exists an element $x \in \mathbb Z_p[\zeta_p]$ such that:
$$Zm=xm.$$
Now set:
$$D(Z)=\prod_{\sigma \in \Delta}(Z-x^{\sigma})\in \Lambda.$$
Then:
$$D(Z)M_F^H=\{ 0\}.$$
Therefore $G(Z)$ divides $D(Z)$ in $\Lambda.$ Thus there exists $\sigma \in \Delta$ such that:
$$x^{\sigma }=j_{1,a}(P)-1.$$
But:
$$C_0(F)^H\simeq \frac{M_F^H}{ZM_F^H}\simeq \frac{\mathbb Z_p[\zeta_p]}{x\mathbb Z_p[\zeta_p]}.$$
Therefore, we get:
$$
J_k(\widetilde{\mathbb F_{\ell}})^{H_a}\otimes_{\mathbb Z} \mathbb Z_p\simeq \frac{\mathbb Z_p[\zeta_p]}{(j_{1,a}(P)-1)\mathbb Z_p[\zeta_p]}.$$
Recall that $j_{1,a}(P)\equiv 1\pmod{\pi^2}.$ Thus:
$$v_p(j_{1,a}(P)-1)=v_p({\rm Log}_p(j_{1,a}(P))).$$
Now:
$${\rm Log}_p(j_{1,a}(P))=\frac{1}{2}f{\rm Log}_p(\ell )+\sum_{\psi \in \widehat{\Delta},\psi \, {\rm odd}}e_{\psi}{\rm Log}_p(j_{1,a}(P)),$$
where $f$ is the order of $\ell$ in $(\mathbb Z/p\mathbb Z)^*.$ Let $\psi \in \widehat \Delta,$ $\psi=\omega^n,$ $n$ odd. If $e_{\psi}{\rm Log}_p(j_{1,a}(P))\not =0,$ then:
$$v_p(e_{\psi}{\rm Log}_p(j_{1,a}(P)))\equiv \frac{n}{p-1} \pmod{\mathbb Z},$$
furthermore:
$$v_p(e_{\psi}{\rm Log}_p(j_{1,a}(P)))>\frac{n}{p-1}\, \Leftrightarrow e_{\psi}j_{1,a}(P)\in U^p.$$
Note also that:
$$v_p(e_{\omega}{\rm Log}_p(j_{1,a}(P)))>\frac{1}{p-1}.$$
The Theorem follows. $\diamondsuit$\par
\newtheorem{Corollary5}[Lemma11]{Corollary}
\begin{Corollary5}\label{Corollary5}
Let $n\in \{ 3,\cdots, p-2\} ,$ $n$ odd.  Let $a \in \{ 1,\cdots ,p-2\}$ such that $1+a^n-(1+a)^n \not \equiv 0\pmod{p}.$ The following assertions are equivalent:\par
\noindent 1) $A(\omega^{1-n})=\{ 0\},$\par
\noindent 2) there exists a prime number $\ell,$ $\ell \not =p,$ such that $\widetilde{d}_p J_k^{H_a}=n.$
\end{Corollary5}
\noindent{\sl Proof} Observe that 2) implies 1) by the above Theorem and Theorem \ref{Theorem3}. Write $\psi =\omega ^n.$ Let $\ell$ be a prime number, $\ell \not =p.$ Write:
$$\mathbb F_{(\ell )}=\frac{\cal{O}}{\ell \cal{O}},$$
and:
$$D_{\ell}=\frac {\mathbb F_{(\ell )}^*}{(\mathbb F_{(\ell )}^*)^p}.$$
Observe that $D_{\ell}$ is a $\mathbb Z_p[\Delta ]$-module. Let $Cyc$ be the group of cyclotomic units of $L.$ We denote the image of $Cyc$ in $D_{\ell}$ by $\overline{Cyc}^{\ell}.$  Then Theorem \ref{Theorem2} asserts that $e_{\psi} \overline{Cyc}^{\ell}=\{ 1\}$ in $D_{\ell }$ if and only if $e_{\psi}j_{1,a}(P)\in U^p,$ where $P$ is a prime of $\cal O$ above $\ell.$ Let:
$$B=L({}^p\sqrt{Cyc}).$$
We assume that 1) holds.  We apply the Chebotarev density theorem to the extension $B/L,$ then there exist infinitely many primes $\ell$ such that:\par
\noindent - $e_{\rho} \overline{Cyc}^{\ell}=\{1 \}$ for $\rho \not = \psi,$\par
\noindent - $e_{\psi} \overline{Cyc}^{\ell}\not =\{1 \}.$\par
\noindent  It remains to apply Theorem \ref{Theorem4} and the above remarks to get 2).  $\diamondsuit$\par
Now, let $\ell$ be a prime number. Let $p$ be an odd  prime number, $p\not =\ell.$ Let $T$ be an indeterminate over $\mathbb F_{\ell}$ and let $E_p/\mathbb F_{\ell} (T)$ be the imaginary quadratic extension defined by:
$$E_p=\mathbb F_{\ell} (T, X)\, {\rm where}\, 
X^2-X+T^p=0.$$
Let $n $ be an odd integer, $n\geq 3.$  Let $S_n (\ell )$ denote the set of primes $p$ such that $\widetilde {d}_pJ_{E_p}=n.$ By our above results, we remark that if $p\in S_n(\ell)$ then  $A(\omega^{1-n})=\{0\}.$ Observe that if $\ell^n\not \equiv 1\pmod{p}$ then $p\not \in S_n(\ell ),$ and therefore $S_n(\ell )$ is a finite set. Set $S(\ell )=\cup_{n}S_{n}(\ell),$ where $n$ runs through the odd integers. Observe that if the order of $\ell$ modulo $p$ is even then $p\not \in S(\ell .)$ Therefore, by a classical result of Hasse (see \cite{LAG}) there exist infinitely many prime $p$ not in $S(\ell )$ (in fact  at least "$2/3$ of the prime numbers "are not in $S(\ell )$). Thus, we ask the following question:\par
\centerline{ is $S(\ell )$ infinite?}\par


 \end{document}